\documentclass[a4paper,12pt]{amsart}
\usepackage{a4wide,graphicx,amssymb}
\usepackage{tikz}
\usetikzlibrary{decorations.pathreplacing,arrows}
\usetikzlibrary{shapes}
 \usepackage{enumerate}
 \usepackage{multirow}
 \usepackage{multicol}

\newcommand{\PP}{\mathbb{P}}

\newcommand{\rbibd}{\textnormal{RBIBD}}
\def\Cr{{\mbox {\sc cr}}}

\def\rcro{{\mbox {\sc rcr}}}

\newif\ifdetails
\detailstrue
\newcommand{\DETAIL}[1]%
{\ifdetails\par\fbox{\begin{minipage}{0.9\linewidth}\textit{Detail:}
      #1\end{minipage}}\par\fi}
\newcommand{\TODO}[1]%
{\ifdetails\par\fbox{\begin{minipage}{0.9\linewidth}\textbf{TODO:}
      #1\end{minipage}}\par\fi}

\newtheorem{lemma}{Lemma}

\newtheorem{theorem}[lemma]{Theorem}

\newcommand{\Tp}{\operatorname{T}}

\newcommand{\old}[1]{{}}

\title{Using block designs in crossing number bounds}

\author{John Asplund} 
\author{\'Eva Czabarka}
\author{Gregory Clark}
\author{Garner Cochran}
\author{Arran Hamm}
\author{Gwen Spencer}
\author{L\'aszl\'o Sz\'ekely}
\author{Libby Taylor}
\author{Zhiyu Wang} 

\address{John Asplund\\Dalton State College\\ Department of Technology and Mathematics\\ 650 College Dr\\ Dalton GA 30720 \\ USA}
\email{jasplund@daltonstate.edu}
\address{\'Eva Czabarka\\Department of Mathematics \\ University of South Carolina \\ Columbia SC 29212 \\ USA
\and Visiting Professor\\ Department of Pure and Applied Mathematics\\ University of Johannesburg\\
P.O. Box 524, Auckland Park, Johannesburg 2006\\South Africa}
\email{czabarka@math.sc.edu}
\address{Gregory Clark\\Department of Mathematics \\ University of South Carolina \\ Columbia SC 29212 \\ USA}
\email{gjclark@math.sc.edu}
\address{Garner Cochran\\Department of Mathematics and Computer Science\\ Berry College\\     2277 Martha Berry Hwy NW\\  Mt Berry GA 30149\\ USA}
\email{gcochran@math.sc.edu}
\address{Arran Hamm\\Winthrop University\\ Department of Mathematics\\ 701 Oakland Ave\\ Rock Hill SC 29733\\ USA}
\email{hamma@winthrop.edu}
\address{Gwen Spencer\\Department of Mathematics and Statistics\\ Smith College\\ Northampton MA 01063 \\ USA}
\email{gspencer@smith.edu}
\address{L\'aszl\'o Sz\'ekely\\Department of Mathematics \\ University of South Carolina \\ Columbia SC 29212 \\ USA
\and Visiting Professor\\ Department of Pure and Applied Mathematics\\ University of Johannesburg\\
P.O. Box 524, Auckland Park, Johannesburg 2006\\ South Africa}
\email{szekely@math.sc.edu}
\address{Libby Taylor\\Department of Mathematics \\ Stanford University\\  Building 380\\ Stanford  CA 94305 \\ USA}
\email{libbyrtaylor@gmail.com}
\address{Zhiyu Wang\\ Department of Mathematics \\ University of South Carolina \\ Columbia SC 29212 \\ USA}
\email{zhiyuw@math.sc.edu}

\subjclass[2010]{Primary 05C10; secondary 05B05, 05D40, 05C85}
\keywords{$k$-planar crossing number, resolvable balanced incomplete block design, Kirkman triple system, resolvable group divisible design}
\thanks{This material is based upon work that started at the Mathematics Research Communities workshop ``Beyond Planarity: Crossing Numbers of Graphs'', Snowbird, Utah, June 11--18, 2017, organized by the American Mathematical Society,  and continued at the 
follow-up workshop ``Biplanar Crossing Numbers and Random Graphs", February 22--25, 2018,
at the University of South Carolina,  with the 
support of the National Science Foundation  contract DMS 1641020 and the American Mathematical Society.
The seventh and the ninth authors were also supported in part by the National Science Foundation  contract  DMS  1600811.}

\begin{document}

\begin{abstract}
The {\em crossing number} $\Cr(G)$ of a graph $G=(V,E)$ is
the smallest number of edge crossings over all drawings of $G$ in the
plane. For any $k\ge 1$, the {\em $k$-planar crossing number} of $G$,
$\Cr_k(G)$, is defined as the minimum of
$\Cr(G_1)+\Cr(G_2)+\ldots+\Cr(G_{k})$ over all graphs $G_1, G_2,\ldots,
G_{k}$ with $\cup_{i=1}^{k}G_i=G$.
Pach et al. [\emph{Computational Geometry: Theory and Applications}
{\bf 68} 2--6, (2018)] showed that for every $k\ge 1$, we have $\Cr_k(G)\le
\left(\frac{2}{k^2}-\frac1{k^3}\right)\Cr(G)$ and that this bound does not remain
true if we replace the constant $\frac{2}{k^2}-\frac1{k^3}$ by any number
smaller than $\frac1{k^2}$. We improve the upper bound to $\frac{1}{k^2}(1+o(1))$ as $k\rightarrow \infty$. For the class of bipartite graphs, we show that 
the best constant is exactly $\frac{1}{k^2}$ for every $k$.
The results extend to the rectilinear
variant of the $k$-planar crossing number.
\end{abstract}

\maketitle

\section{Introduction}

This note improves on results of Pach, Sz\'ekely, T\'oth, and T\'oth \cite{pach}.  We follow the introduction
of that paper.

A \emph{drawing} of a graph $G=(V,E)$ is a planar representation of $G$ such
that every vertex $v\in V$ corresponds to a point of the plane and every edge
$uv\in E$ is represented by a simple continuous curve between the points
corresponding to $u$ and $v$, which does not pass through any point
representing a vertex of $G$. We assume for simplicity
that no two curves share infinitely many points, no two curves are
tangent to each other, and no three curves pass through the same
point. The \emph{crossing number} $\Cr(G)$ of $G$ is defined as the minimum number of
edge crossings in a drawing of $G$. For surveys,
see~\cite{schaefer,success}, and the recent monograph \cite{schaefer2}. 
Clearly, $G$ is planar if and only if $\Cr(G)=0$.

Selfridge (see~\cite{Ha61}) noticed that by Euler's polyhedral formula,
$K_{11}$, the complete graph on $11$ vertices, cannot be written as the union
of two planar graphs. Battle, Harary, and Kodama~\cite{BaHK62} and
independently Tutte~\cite{Tu63a} proved that the same is true for $K_9$, but
not for $K_8$. This led Tutte~\cite{Tu63b} to introduce the
\emph{thickness} of a graph $G$, which is the minimum number of planar graphs
that $G$ can be decomposed into. The notion is relevant for VLSI
chip design, where it corresponds to the number of layers required for
realizing a network so that there is no crossing within a layer (see
Mutzel, Odenthal, and Scharbrodt~\cite{MuOS98} for a survey). If the thickness
of $G$ is at most $2$, $G$ is called \emph{biplanar}. Mansfield proved that it
is an NP-complete problem to decide whether a graph is biplanar; see~\cite{beineke,Ma83}.

Owens~\cite{owens} defined the \emph{biplanar crossing number} $\Cr_2(G)$ of $G$ 
 as the minimum sum of the crossing numbers of two graphs,
$G_0$ and $G_1$, whose union is $G$. $G$ is
biplanar precisely when its biplanar crossing number is $0$.
Shahrokhi et al.~\cite{kplanar} extended this notion as follows. For any positive integer $k\ge 1$, define the
\emph{$k$-planar crossing number} $\Cr_k(G)$ of $G$ as the minimum of
$\Cr(G_1)+\Cr(G_2)+\ldots+\Cr(G_{k})$,  where the minimum is taken over
all graphs $G_1, G_2,\ldots, G_{k}$ whose union is $G$, that is,
$\bigcup_{i=1}^{k}E(G_i)=E(G)$.

Spencer~\cite{spencer} showed that for sufficiently large $c$ for all $p>c/n$ with high probability
the biplanar crossing number of Erd\H{o}s-R\'enyi
random graphs  is
$\Theta(n^4p^2)$, and claimed a similar result for $k$-planar crossing 
numbers without proof. Asplund et al.~\cite{kpl} gave a proof and extended this result for random $d$-regular graphs, where $d$ exceeds a certain 
threshold.

Czabarka, S\'ykora,
Sz\'ekely, and Vr\v to  \cite{bipII} proved that for every
graph $G$, we have
\begin{equation} \label{3/8}
\Cr_2(G)\leq \frac{3}{8}\Cr(G).
\end{equation}
They also showed~\cite{bipsurvey} that this inequality does not remain true
if the constant $\frac{3}{8}=0.375$ is replaced by anything less than $\frac{8}{119}\approx 0.067$.

Pach et al.~\cite{pach}  extended this investigation to the relationship between the $k$-planar
crossing number and the (ordinary) crossing number of a graph.
For every integer  $k\ge 1$, they defined
$$\alpha_k= \sup \frac{\Cr_k(G)}{\Cr(G)},$$
where the supremum is taken over all \emph{nonplanar} graphs $G$. 
The results mentioned from~\cite{bipII} yield $0.067<\alpha_2\le \frac38=0.375$. Pach et al.~\cite{pach} proved that
for every positive integer $k$, 
\begin{equation}\label{eq:prevbest} \frac{1}{k^2}\le\alpha_k\le \frac{2}{k^2}-\frac{1}{k^3}.\end{equation}
Note that for $k=2$, (\ref{eq:prevbest}) returns the value 3/8 given in (\ref{3/8}), and the present paper does not improve this upper bound on $\alpha_2$ either. In this paper, we show that the lower bound in (\ref{eq:prevbest}) is asymptotically correct as $k\rightarrow \infty$. 
\begin{theorem}\label{th:main}
$\alpha_k=\frac{1}{k^2}(1+o(1))$ as $k\rightarrow \infty$.\end{theorem}

As Theorem~\ref{th:main} and its proof surrender control over the $o(1)$ term, it is of interest to determine the values of $\alpha_k$ for small $k$. 
To this end,
we improve the upper bound $\frac{2}{k^2}-\frac{1}{k^3}$ for  $3\leq k\leq 10$, see Table~\ref{tabla2}, using the following theorem.
\begin{theorem}\label{th:main2} We have
\begin{enumerate}[{\rm (i)}]
\item\label{th:4}  $\alpha_4\leq \frac{235}{2401}$;
\item\label{th:2mod3}       $\alpha_k\leq \frac{12k-11}{(2k-1)^3} $ for $k\equiv 2 \pmod{3}$;
\item\label{th:2mod4}       $\alpha_k\leq \frac{36k-35}{(3k-2)^3} $ for $k\equiv 2 \pmod{4}$;
\item\label{th:0mod3}   $\alpha_k\leq \frac{3k-1}{2k^3}$ for $k\equiv 0 \pmod{3}$;
\item\label{th:odd}  $\alpha_k\leq \frac{2}{k(k+1)}$ for odd $k$.
\end{enumerate}
\end{theorem}
Note that while for odd $k$  the expression in (\ref{th:odd}) offers an improvement over (\ref{eq:prevbest}) that is in diminishing proportion as 
$k\rightarrow \infty$, the gain is still meaningful for small values of $k$. 
In contrast, 
(\ref{th:2mod3}), (\ref{th:2mod4}), and (\ref{th:0mod3}) also offer an asymptotic improvement over (\ref{eq:prevbest}). 

\begin{table} \label{tabla2}
\begin{center}
\def\arraystretch{2}
\begin{tabular}{|l|c|c|c|}
\hline
$k$ & $\alpha_k$ bound from ~\eqref{eq:prevbest} & $\alpha_k$ bound improved & lower bound\\
\hline\hline
$3$ &$\frac{5}{27}\lessapprox 0.1852$&$\frac{1}{6}\lessapprox 0.1667$ (\ref{th:odd})&$\frac{1}{9} \gtrapprox 0.1111$\\
\hline
$4$ &$\frac{7}{64}\lessapprox 0.1094$&$\frac{235}{2401}\lessapprox 0.0979$ (\ref{th:4})&$\frac{1}{16}=0.0625$\\
\hline
$5$ &$\frac{9}{125}= 0.072$&$\frac{1}{15}\lessapprox 0.0667$ (\ref{th:odd})&$\frac{1}{25} =0.04 $\\
\hline
$6$ &$\frac{11}{216}\lessapprox 0.0510$&$\frac{17}{432}\lessapprox 0.0394$ (\ref{th:0mod3})&$\frac{1}{36} \gtrapprox 0.0277$\\
\hline
$7$ &$\frac{13}{343}\lessapprox 0.0380$&$\frac{1}{28}\lessapprox 0.0358$ (\ref{th:odd})&$\frac{1}{49} \gtrapprox 0.2040$\\
\hline
$8$ &$\frac{15}{512}\lessapprox 0.0293$&$\frac{85}{3375}\lessapprox 0.0252$ (\ref{th:2mod3})&$\frac{1}{64}=0.015625$\\
\hline
$9$ &$\frac{17}{729}\lessapprox 0.0234$&$\frac{13}{729}\lessapprox 0.0179$ (\ref{th:0mod3})&$\frac{1}{81} \gtrapprox 0.0123 $\\
\hline
$10$ &$\frac{19}{1000}=0.019$&$\frac{325}{21952}\lessapprox 0.0149$ (\ref{th:2mod4})&$\frac{1}{100} = 0.01$\\
\hline
\end{tabular}
\end{center}
\caption{Comparison of the best upper bounds for $\alpha_k$ from \eqref{eq:prevbest}, due to Pach et al. \cite{pach},   our upper bounds, and the lower bound $\frac{1}{k^2}$ for $3\leq k\leq 10$. Roman numerals in the second column refer to cases of Theorem~\ref{th:main2}.}–
\end{table}

We also consider the restriction of the problem to bipartite graphs.  To this end, define
$$\beta_k= \sup \frac{\Cr_k(G)}{\Cr(G)},$$
where the supremum is taken over all \emph{nonplanar bipartite} graphs $G$. When restricted to bipartite graphs, we can show that the lower bound in (\ref{eq:prevbest}) is exact.

\begin{theorem}\label{th:mainbip}
For all $k$, 
$\beta_k=\frac{1}{k^2}$.
\end{theorem}

The \emph{rectilinear crossing number}, $\rcro(G)$, of a graph $G$ is the minimum number of crossings
over all \emph{straight-line} drawings of $G$, in which the edges are represented by line segments. Obviously,
we have $\Cr(G)\leq \rcro(G)$ for every graph $G$. For every $t\geq 4$, Bienstock and Dean~\cite{BD93} constructed families of graphs whose crossing number is at most $t$ and whose rectilinear crossing number is unbounded.

Similar to $\Cr_k(G)$, we define the \emph{rectilinear $k$-planar crossing number} of a graph $G$, denoted $\rcro_k(G)$, as the minimum of $\rcro(G_1)+\rcro(G_2)+\ldots+\rcro(G_{k})$, where the minimum is taken over
all graphs $G_1, G_2,\ldots, G_{k}$ whose union is $G$. It is likewise clear that $\Cr_k(G)\leq \rcro_k(G)$ for
every  positive integer $k$.  The analogue of $\alpha_k$ 
 is
$$\overline{\alpha}_k= \sup \frac{\rcro_k(G)}{\rcro(G)},$$
where the supremum is taken over all \emph{nonplanar} graphs $G$, and the analogue of $\beta_k$ 
is
$$\overline{\beta}_k= \sup \frac{\rcro_k(G)}{\rcro(G)},$$
where the supremum is taken over all \emph{bipartite nonplanar} graphs $G$  (as planar is the same as rectilinear planar by \cite{f}). We have 
\begin{theorem}\label{th:main3} Theorems~\ref{th:main} and~\ref{th:main2} remain true if we replace $\alpha_k$ with $\overline{\alpha}_k$ and
$\beta_k$ with $\overline{\beta}_k$, consequently the bounds for $\alpha_k$ in Table~\ref{tabla2} apply for $\overline{\alpha}_k$ as well.
\end{theorem}

\section{Methodology} \label{method}
We generalize the procedure that was defined for two planes in \cite{bipII} and extended to $k$ planes in \cite{pach}.
 Given an integer $k>1$, we create a $k$-planar drawing of $G$ in the following way.
We number the $k$ planes with $1,2,\ldots,k$, and describe a probabilistic procedure that assigns a plane to each edge, resulting in a graph $G_i$ on the $i$-th plane.
Let $K_s^{o}$ denote the complete graph on the vertex set $\{1,2,...,s\}$, with a loop edge added at every vertex and let the graphs $H_1,H_2,...,H_k$ partition the edge set of $K_s^o$. (In \cite{bipII}, where $k=2$, the choice
of $s$ was  2, with $H_1=$ a single edge and $H_2$ = two loops; in \cite{pach} the choice was $s=k$, for odd $k$  every $H_i$ consisted 
of a single loop and a perfect matching on the remaining $k-1$ vertices, while for even $k$  every $H_i$ was either
a perfect matching, or two loops and a perfect matching on the remaining $k-2$ vertices.) We call the components of $H_1,H_2,...,H_k$ {\it types}. 
Note that two distinct components (no matter whether they are in the same $H_i$ or not) are different types even if they are isomorphic. For an example, see Figure \ref{fig:4plane}.
We define the {\it type} of an edge of $K_s^{o}$
(either loop or not) as  the unique component of the subgraph $H_i$ containing it.

Given a graph $G$, we distribute the edges of $G$ into the $k$ planes as follows. To each vertex $v$ of $G$ assign a value $\xi(v)$ randomly and uniformly chosen from $\{1,2,...,s\}$ where $s$ is chosen carefully depending on values of $k$.
If $uv$ is an edge of $G$, assign the $uv$  edge to the $j$-th plane if $\{\xi(u),\xi(v)\}\in E(H_j)$ (where $\{i\}=\{i,i\}$ is the loop on vertex $i$ of $K_s^{o}$). 
As the $E(H_j)$'s partition $E(K_s^{o})$, there is exactly one $j$ assigned to each edge.  

We will use an optimal drawing $\mathcal{D}$ of $G$ realizing $\Cr(G)$ to create a $k$-planar drawing. 
It is well-known that in $\mathcal{D}$  every pair $e,f$ of crossing edges has four distinct endvertices and the edges $e,f$ have exactly one point in common
\cite{schaefer, schaefer2, success}. Denote by $G_j$ the subgraph of $G$ containing the edges assigned to the $j$-th plane. Draw $G_j$ in the $j$-th plane
following the drawing $\mathcal{D}$, i.e. the drawing of each edge $uv$ in $G_j$ follows the curve representing the $uv$ edge in $\mathcal{D}$.

Assume that $C$ is a component of $H_i$. Then clearly the subgraph of $G$ induced by $\{v\in V(G): \xi(v)\in V(C)\}$ is a union of components in $G_i$. We modify
our $k$-planar drawing
to further reduce the number of crossing edge pairs by translating the drawings of subgraphs of $G_i$ on the vertex sets $\{v\in V(G): \xi(v)\in V(C)\}$ for the components of $H_i$ far enough from each other so that if $e_1,e_2\in E(G_i)$ and vertices of $e_1$ and $e_2$ are mapped to different components of $H_i$ by $\xi$,
then curves corresponding to $e_1$ and $e_2$ do not cross in
the drawing of $G_i$. 

Assume that $uv$ and $wz$ are a pair of crossing edges in the optimal drawing $\mathcal{D}$ of $G$, and hence have 4 distinct
endpoints. The probability that this edge pair is still crossing in the random $k$-planar drawing above, is exactly 
\begin{equation} \label{type}
q=\PP[{\rm type}(\xi(u),\xi(v))={\rm type}(\xi(w),\xi(z))].
\end{equation}
The value of $q$ does not depend on which crossing edge pair $uv$ and $wz$ was selected from $\mathcal{D}$, so 
the expected number of crossings in our random $k$-planar drawing is  
\begin{equation} \label{use}
q\Cr(G).
\end{equation}
It follows that there exists a $k$-planar drawing of $G$ which has at most  $q\Cr(G)$ crossings. If this holds for a particular $q$ for all graphs $G$,
then we establish
\begin{equation} \label{alphabound}
\alpha_k\leq q.
\end{equation}
Note that this method can be further enhanced by replacing the base graph $K_s^{o}$ by a graph that is missing some edges (but modifying $\xi$ so no edge
of $G$ is matched to a missing edge) and allowing an edge of the base graph to appear in several of the $H_i$ (and employing another probabilistic procedure to decide
which plane we assign an edge $uv$ to when $\xi(u)\xi(v)$ appears in several of the $H_i$).
We will make use of these modifications in Sections~\ref{odd} and \ref{bipartite}, where we discuss them in more detail. For a warm-up,
we start with $k=4$.

\section{Proof to Theorem~\ref{th:main2}({\normalfont\ref{th:4}}): the case  $k=4$} \label{four}
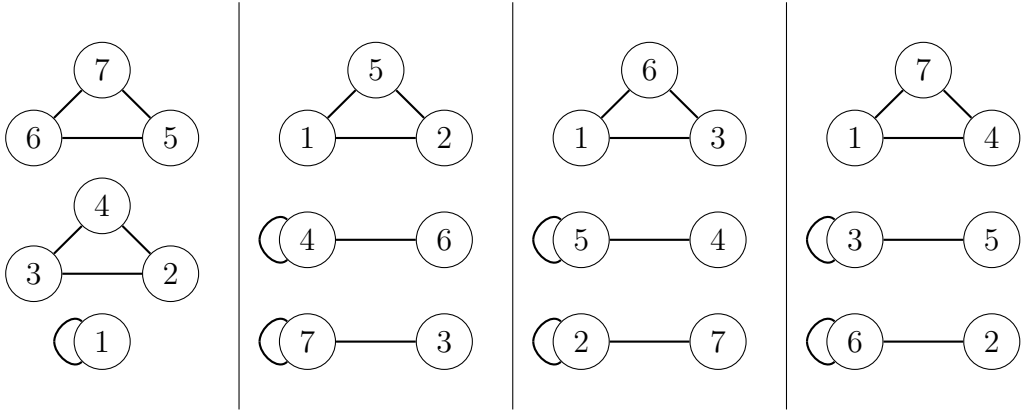
\begin{figure}[htbp]
\begin{center}
\begin{tikzpicture}[ 
    pile/.style={thick, -,black, =stealth'},scale=.9]    
        \node(a1)[shape=circle,draw, fill=white, inner sep=4pt]  at (0,0) {$1$}; 
        \node(b1)[shape=circle,draw, fill=white,  inner sep=4pt]  at (1,1) {$2$};
        \node(c1)[shape=circle,draw, fill=white,  inner sep=4pt]  at (-1,1) {$3$};
        \node(d1)[shape=circle,draw, fill=white,  inner sep=4pt]  at (0,2) {$4$}; 
        \node(e1)[shape=circle,draw, fill=white,  inner sep=4pt]  at (1,3) {$5$}; 
        \node(f1)[shape=circle,draw, fill=white,  inner sep=4pt]  at  (-1,3) {$6$};
        \node(g1)[shape=circle,draw, fill=white,  inner sep=4pt]  at (0,4) {$7$};
        \draw[pile] (b1)--(c1)--(d1)--(b1);
        \draw[pile] (e1)--(f1)--(g1)--(e1);
       \path[pile,out=135,in=90] (a1) edge (-.7,0);     
       \path[pile,out=-135,in=-90] (a1) edge (-.7,0);                
        \draw (2,-1)--(2,5); %
        \node(a2)[shape=circle,draw, fill=white,  inner sep=4pt]  at (3,3) {$1$}; 
        \node(b2)[shape=circle,draw, fill=white,  inner sep=4pt]  at (5,3) {$2$};
        \node(c2)[shape=circle,draw, fill=white,  inner sep=4pt]  at (5,0) {$3$};
        \node(d2)[shape=circle,draw, fill=white,  inner sep=4pt]  at (3,1.5) {$4$};
        \node(e2)[shape=circle,draw, fill=white,  inner sep=4pt]  at (4,4) {$5$};
        \node(f2)[shape=circle,draw, fill=white,  inner sep=4pt]  at  (3,0) {$7$};
        \node(g2)[shape=circle,draw, fill=white,  inner sep=4pt]  at (5,1.5) {$6$};
       \draw[pile] (b2)--(e2)--(a2)--(b2);
        \draw[pile] (d2)--(g2);
        \draw[pile] (c2)--(f2);   
        \path[pile,out=135,in=90] (d2) edge (2.3,1.5);     
       \path[pile,out=-135,in=-90] (d2) edge (2.3,1.5);                
       \path[pile,out=135,in=90] (f2) edge (2.3,0);     
       \path[pile,out=-135,in=-90] (f2) edge (2.3,0);                                     
         \draw (6,-1)--(6,5); %
        \node(a3)[shape=circle,draw, fill=white,  inner sep=4pt]  at (7,3) {$1$}; 
        \node(b3)[shape=circle,draw, fill=white,  inner sep=4pt]  at (9,0) {$7$};
        \node(c3)[shape=circle,draw, fill=white,  inner sep=4pt]  at (9,3) {$3$};
        \node(d3)[shape=circle,draw, fill=white,  inner sep=4pt]  at (9,1.5) {$4$}; 
        \node(e3)[shape=circle,draw, fill=white,  inner sep=4pt]  at (7,1.5) {$5$};
        \node(f3)[shape=circle,draw, fill=white,  inner sep=4pt]  at  (8,4) {$6$};
        \node(g3)[shape=circle,draw, fill=white,  inner sep=4pt]  at (7,0) {$2$};
       \draw[pile] (c3)--(f3)--(a3)--(c3);
        \draw[pile] (b3)--(g3);
        \draw[pile] (d3)--(e3);       
        \path[pile,out=135,in=90] (e3) edge (6.3,1.5);     
       \path[pile,out=-135,in=-90] (e3) edge (6.3,1.5);                
       \path[pile,out=135,in=90] (g3) edge (6.3,0);     
       \path[pile,out=-135,in=-90] (g3) edge (6.3,0);                                                      
        \draw (10,-1)--(10,5); %
        \node(a4)[shape=circle,draw, fill=white,  inner sep=4pt]  at (11,3) {$1$}; 
        \node(b4)[shape=circle,draw, fill=white,  inner sep=4pt]  at (11,0) {$6$};
        \node(c4)[shape=circle,draw, fill=white,  inner sep=4pt]  at (11,1.5) {$3$};
        \node(d4)[shape=circle,draw, fill=white,  inner sep=4pt]  at (13,3) {$4$}; 
        \node(e4)[shape=circle,draw, fill=white,  inner sep=4pt]  at (13,1.5) {$5$};
        \node(f4)[shape=circle,draw, fill=white,  inner sep=4pt]  at  (13,0) {$2$};
        \node(g4)[shape=circle,draw, fill=white,  inner sep=4pt]  at (12,4) {$7$};
       \draw[pile] (d4)--(g4)--(a4)--(d4);
        \draw[pile] (b4)--(f4);
        \draw[pile] (c4)--(e4);                
        \path[pile,out=135,in=90] (c4) edge (10.3,1.5);     
       \path[pile,out=-135,in=-90] (c4) edge (10.3,1.5);                
       \path[pile,out=135,in=90] (b4) edge (10.3,0);     
       \path[pile,out=-135,in=-90] (b4) edge (10.3,0);                                                      

\end{tikzpicture}
\old{\begin{tikzpicture}[ 
    pile/.style={thick, -,black, =stealth'}]    
        \node(a1)[shape=circle,draw, fill=white, inner sep=4pt]  at (0,0) {$0$}; 
        \node(b1)[shape=circle,draw, fill=white,  inner sep=4pt]  at (1,1) {$1$};
        \node(c1)[shape=circle,draw, fill=white,  inner sep=4pt]  at (-1,1) {$2$};
        \node(d1)[shape=circle,draw, fill=white,  inner sep=4pt]  at (0,2) {$3$}; 
        \node(e1)[shape=circle,draw, fill=white,  inner sep=4pt]  at (1,3) {$4$}; 
        \node(f1)[shape=circle,draw, fill=white,  inner sep=4pt]  at  (-1,3) {$5$};
        \node(g1)[shape=circle,draw, fill=white,  inner sep=4pt]  at (0,4) {$6$};
        \draw[pile] (b1)--(c1)--(d1)--(b1);
        \draw[pile] (e1)--(f1)--(g1)--(e1);
       \path[pile,out=135,in=90] (a1) edge (-.7,0);     
       \path[pile,out=-135,in=-90] (a1) edge (-.7,0);                
        \draw (2,-1)--(2,5); %
        \node(a2)[shape=circle,draw, fill=white,  inner sep=4pt]  at (3,3) {$0$}; 
        \node(b2)[shape=circle,draw, fill=white,  inner sep=4pt]  at (5,3) {$1$};
        \node(c2)[shape=circle,draw, fill=white,  inner sep=4pt]  at (5,0) {$2$};
        \node(d2)[shape=circle,draw, fill=white,  inner sep=4pt]  at (3,1.5) {$3$};
        \node(e2)[shape=circle,draw, fill=white,  inner sep=4pt]  at (4,4) {$4$};
        \node(f2)[shape=circle,draw, fill=white,  inner sep=4pt]  at  (3,0) {$5$};
        \node(g2)[shape=circle,draw, fill=white,  inner sep=4pt]  at (5,1.5) {$6$};
       \draw[pile] (b2)--(e2)--(a2)--(b2);
        \draw[pile] (d2)--(g2);
        \draw[pile] (c2)--(f2);   
        \path[pile,out=135,in=90] (d2) edge (2.3,1.5);     
       \path[pile,out=-135,in=-90] (d2) edge (2.3,1.5);                
       \path[pile,out=135,in=90] (f2) edge (2.3,0);     
       \path[pile,out=-135,in=-90] (f2) edge (2.3,0);                                     
         \draw (6,-1)--(6,5); %
        \node(a3)[shape=circle,draw, fill=white,  inner sep=4pt]  at (7,3) {$0$}; 
        \node(b3)[shape=circle,draw, fill=white,  inner sep=4pt]  at (9,0) {$1$};
        \node(c3)[shape=circle,draw, fill=white,  inner sep=4pt]  at (9,3) {$2$};
        \node(d3)[shape=circle,draw, fill=white,  inner sep=4pt]  at (9,1.5) {$3$}; 
        \node(e3)[shape=circle,draw, fill=white,  inner sep=4pt]  at (7,1.5) {$4$};
        \node(f3)[shape=circle,draw, fill=white,  inner sep=4pt]  at  (8,4) {$5$};
        \node(g3)[shape=circle,draw, fill=white,  inner sep=4pt]  at (7,0) {$6$};
       \draw[pile] (c3)--(f3)--(a3)--(c3);
        \draw[pile] (b3)--(g3);
        \draw[pile] (d3)--(e3);       
        \path[pile,out=135,in=90] (e3) edge (6.3,1.5);     
       \path[pile,out=-135,in=-90] (e3) edge (6.3,1.5);                
       \path[pile,out=135,in=90] (g3) edge (6.3,0);     
       \path[pile,out=-135,in=-90] (g3) edge (6.3,0);                                                      
        \draw (10,-1)--(10,5); %
        \node(a4)[shape=circle,draw, fill=white,  inner sep=4pt]  at (11,3) {$0$}; 
        \node(b4)[shape=circle,draw, fill=white,  inner sep=4pt]  at (11,0) {$1$};
        \node(c4)[shape=circle,draw, fill=white,  inner sep=4pt]  at (11,1.5) {$2$};
        \node(d4)[shape=circle,draw, fill=white,  inner sep=4pt]  at (13,3) {$3$}; 
        \node(e4)[shape=circle,draw, fill=white,  inner sep=4pt]  at (13,1.5) {$4$};
        \node(f4)[shape=circle,draw, fill=white,  inner sep=4pt]  at  (13,0) {$5$};
        \node(g4)[shape=circle,draw, fill=white,  inner sep=4pt]  at (12,4) {$7$};
       \draw[pile] (d4)--(g4)--(a4)--(d4);
        \draw[pile] (b4)--(f4);
        \draw[pile] (c4)--(e4);                
        \path[pile,out=135,in=90] (c4) edge (10.3,1.5);     
       \path[pile,out=-135,in=-90] (c4) edge (10.3,1.5);                
       \path[pile,out=135,in=90] (b4) edge (10.3,0);     
       \path[pile,out=-135,in=-90] (b4) edge (10.3,0);                                                      

\end{tikzpicture}
}
\end{center}
\caption{Partitioning $K_7^o$ into 12 types  in $4$ planes.
}  \label{fig:4plane}
\end{figure}

Choose $s=7$, 
and see Figure~\ref{fig:4plane} for the partition of $K_7^o$ into 12 types on 4 planes.

\old{
If $|\Tp(uv)|=1$ for some $i$, then assign the edge to plane $0$ when $\xi_u=0$, plane $1$ when $\xi_u\in\{3,5\}$, plane $2$ when $\xi_u\in\{4,6\}$ and plane $3$
when $\xi_u\in\{1,2\}$.
If $|\Tp(uv)|=2$, then 
\begin{enumerate}
\item Assign the edge to plane $0$ if $\xi_u,\xi_v\in\{1,2,3\}$ or $\xi_u,\xi_v\in\{4,5,6\}$.
\item Assign the to plane $1$ if $\xi_u\xi_v\in\{0,1,4\}$ or $\Tp(uv)\in\{\{3,6\},\{2,5\}\}$
\item Assign the edge to plane $2$ if $\xi_u\xi_v\in\{0,2,5\}$ or $\Tp(uv)\in\{\{3,4\},\{1,6\}\}$
\item Assign the edge to plane $3$ in all other cases, i.e. when $\xi_u\xi_v\in\{0,2,6\}$ or $\Tp(uv)\in\{\{2,4\},\{1,5\}\}$
\end{enumerate}
}

Take a crossing edge pair $\{ab,cd\}$. Given $V(G)=\{1,2,...,n\}$, without loss of generality we assume that $a<b,c<d$ and $a<c$. This determines the $(a,b,c,d)$ quadruple uniquely for each crossing edge pair.
We compute the probability that $ab$ and $cd$ still cross in the $k$-planar drawing we provided
by counting the number of ways the edge pair can be labeled and remain crossing, and dividing it by $7^4$, the total number of possible labelings.
\begin{enumerate}
\item{If $\{\xi(a),\xi(b)\}\cap\{\xi(c),\xi(d)\}=\emptyset$, then the edge pair does not remain crossing.}
\item{If $\xi(a)=\xi(b)=\xi(c)=\xi(d)$, then the edge pair remains crossing, and  $7$ different labelings yield such a situation.}
\item{If for some $i\ne j$, $\{\{\xi(a),\xi(b)\},\{\xi(c),\xi(d)\}\}=\{\{i\},\{i,j\}\}$, then the edge pair remains crossing only if $i\ne 1$.
There are  $6\cdot 4=24$ ways to label the vertices this way.}
\item{If for some $i\ne j$, $\{\{\xi(a),\xi(b)\},\{\xi(c),\xi(d)\}\}=\{\{i,j\},\{i,j\}\}$, then the edge pair remains crossing,
and there are $\binom{7}{2}\cdot 4=84$ ways to label the vertices this way.}
\item{If   $\{\{\xi(a),\xi(b)\},\{\xi(c),\xi(d)\}\}=\{\{i,j\},\{i,h\}\}$ for three different  numbers $i,j,h$, then the edge pair remains crossing
when $i,j,h$ appear in some triangle in one of the planes in Figure~\ref{fig:4plane}, which can happen
in $15\cdot 8=120$ ways.}
\end{enumerate}

Summing over all possible outcomes yields the probability that a crossing edge pair remains crossing in this $k$-planar drawing is 
\begin{eqnarray*}
\frac{120+84+24+7}{7^4}=\frac{235}{2401}\lessapprox0.0979. 
\end{eqnarray*}

\section{Resolvable BIBDs and proof of Theorem~\ref{th:main} and~\ref{th:main2}}

A {\it resolvable} BIBD, denoted as $\rbibd(s,\ell,\lambda)$, is a collection $P_1,\ldots,P_m$ of partitions of an underlying $s$-element set
into  $\ell$-element subsets such that
 every 2-element subset of the $s$-element set is contained by exactly $\lambda$ of the  $ms/\ell$ $\ell$-element sets listed in the partitions. We restrict ourselves to $\lambda=1$, that is, each $2$-element subset of the $s$-element set is contained in precisely one of the $\ell$-element sets listed in the partitions. 
 
 Note that the existence of such a design
 implies
 that $|P_i|=\frac{s}{\ell}$ and $m\frac{s}{\ell}\binom{\ell}{2}=\binom{s}{2}$, i.e. $m=\frac{s-1}{\ell-1}$, which gives the
  well known necessary condition that $s\equiv \ell \pmod{\ell(\ell-1)}$ for the existence of such a resolvable BIBD. For the $\ell=2$ case, 
which is the factorization of complete graphs into matchings, this condition is also sufficient. For the $\ell=3$ case (known as a Kirkman triple system) it is also a sufficient condition \cite{kirk}, and for $\ell=4$ the corresponding $s\equiv 4 \pmod{12}$ it is also a sufficient condition \cite{hanani}. For every $\ell$, the congruence is also a sufficient condition for all
  $s>s_0(\ell)$ \cite{RV}. Further,
for every even $\ell\geq 4$, the congruence implies existence for $s>\exp\{\exp\{\ell^{18\ell^2}\}\}$ \cite{chang}. 

Assuming that a $\rbibd(s,\ell,1)$ exists, let $k=m+1$, and 
for $i=1,2,...,m$, let $H_i$ be a disjoint union of $K_\ell$'s, whose vertex sets are the $\ell$-element
sets in the partition classes of the partition $P_i$.  For $i=m+1$, we put the $s$ loops into $H_{m+1}$. Following the drawing argument
in Section~\ref{method}, we evaluate the value of $q$.

Consider a crossing edge pair $\{ab,cd\}$ in $G$ as we did in Section~\ref{four}. The following $\xi$-assignments will leave the edge pair crossing:
\begin{enumerate}
\item $\xi(a)=\xi(b)=\xi(c)=\xi(d)$:  $s$ different labelings of these 4 vertices yield such a situation.

\item For some $i\ne j$, $\{\{\xi(a),\xi(b)\},\{\xi(c),\xi(d)\}\}=\{\{i,j\},\{i,j\}\}$:
there are $ \frac{s-1}{\ell-1}\cdot \frac{s}{\ell}\cdot \binom{\ell}{2}\cdot 4$ ways to label the vertices this way.

\item $\{\{\xi(a),\xi(b)\},\{\xi(c),\xi(d)\}\}=\{\{i,j\},\{i,h\}\}$ for three different  numbers $i,j,h$ that 
appear  together in some $\ell$-set of some partition: there are 
$ \frac{s-1}{\ell-1}\cdot \frac{s}{\ell}\cdot \ell\binom{\ell-1}{2}\cdot 8   $ ways to label the vertices this way.

\item $\{\{\xi(a),\xi(b)\},\{\xi(c),\xi(d)\}\}=\{\{i,j\},\{h,g\}\}$ for four different numbers $i,j,h,g$ that
appear  together in some $\ell$-set of some partition: there are 
$ \frac{s-1}{\ell-1}\cdot \frac{s}{\ell}\cdot 3\binom{\ell}{4}\cdot 8   $ ways to label the vertices this way.
\end{enumerate}
Summing over all possibilities and dividing by $s^4$, the total number of labelings of the four vertices, we obtain
\begin{equation} \label{general}
q=\frac{1+(s-1)(\ell^2-\ell)}{s^3}.
\end{equation}

Now we are ready to show the main theorems. 

\begin{subsection}{Proof of Theorem~\ref{th:main}: $\alpha_k=\frac{1}{k^2}(1+o(1))$ as $k\rightarrow \infty$}

\begin{proof}
Note that $\frac{1}{k^2}\le \alpha_k$ from (\ref{eq:prevbest}), so we only have 
to provide an upper bound. 
We will show that for any $\ell\ge 2$ we have $\alpha_k\le\frac{\ell+1}{\ell-1}\cdot\frac{1}{k^2}(1+o(1))$ as $k\rightarrow\infty$. Letting $\ell\rightarrow\infty$ proves the claim.
To this end, fix an $\ell\ge 2$. For a given $k$, set $s=s_k=(k-1)(\ell-1)+1$  (so $k=\frac{s-1}{\ell-1}+1$).
If an $\rbibd(s,\ell,1)$ exists, then  (\ref{general}) gives
\begin{equation*}
\alpha_k\leq q< \frac{k}{(k-1)^3}\cdot \frac{\ell}{\ell-1}+\frac{1}{(k-1)^3(\ell-1)^3}
< \frac{k}{(k-1)^3}\cdot \frac{\ell+1}{\ell-1}.
\end{equation*}

While this may not be true, we know that if $s'$ is sufficiently large and $s'\equiv\ell\pmod{\ell(\ell-1)}$ then an $\rbibd(s',\ell,1)$ does exist.
This means that for $k$ sufficiently large, there exists an $s'$ such that $s_k\geq s'> s_k-\ell(\ell-1)$ and an
$\rbibd(s',\ell, 1)$ exists. Set $k'=\frac{s'-1}{\ell-1}+1$, an integer (so $s'=s_{k'}$). It is easy to see that $k\geq k'> k-\ell$, and
$$\alpha_k\leq \alpha_{k'}\leq \frac{k'}{(k'-1)^3}\cdot \frac{\ell+1}{\ell-1}\leq \frac{k}{(k-\ell-1)^3}\cdot \frac{\ell+1}{\ell-1}=\frac{\ell+1}{\ell-1} \cdot\frac{1}{k^2}\cdot (1+o(1)),$$
verifying our claim.
\end{proof}
\end{subsection}

Now we turn to the proof of Theorem~\ref{th:main2}. Note that Theorem ~\ref{th:main2}(\ref{th:4}) is already shown in Section ~\ref{four}. 

\begin{subsection}{Proof of Theorem ~\ref{th:main2}(\ref{th:2mod3}): $\alpha_k\leq \frac{12k-11}{(2k-1)^3} $ for $k\equiv 2 \pmod{3}$}
\begin{proof}
In Equation \eqref{general}, choose $\ell=3$, and assume that $k\equiv 2\pmod{3}$. Then 
$k\equiv 2$ or $5\pmod{6}$. Set $s=2k-1$. Easy calculation show that
$s\equiv 3\pmod{6}$, and therefore a Kirkman triplet system exists on $s$ vertices. Equation (\ref{general}) gives
$q=\frac{12k-11}{(2k-1)^3}$, proving Theorem~\ref{th:main2}(\ref{th:2mod3}).
\end{proof}
\end{subsection}

\begin{subsection}{Proof of Theorem ~\ref{th:main2}(\ref{th:2mod4}): $\alpha_k\leq \frac{36k-35}{(3k-2)^3} $ for $k\equiv 2 \pmod{4}$}

\begin{proof}
In Equation \eqref{general}, choose $\ell=4$, and assume that $k\equiv 2\pmod{4}$. Then  
$3k\equiv 6\pmod{12}$. Set $s=3k-2$, giving $s\equiv 4\pmod{12}$, and therefore a resolvable BIBD exists
with $\ell=4$ on $s$ vertices. Equation (\ref{general}) yields $q=\frac{36k-35}{(3k-2)^3}$, proving Theorem~\ref{th:main2}(\ref{th:2mod4}).
\end{proof}
\end{subsection}

\begin{subsection}{Proof of Theorem~\ref{th:main2}(\ref{th:0mod3}): $\alpha_k\leq \frac{3k-1}{{2k}^3} $ for $k\equiv 0 \pmod{3}$}

\begin{proof}

When $k\equiv 0\pmod{3}$, set $s=2k$, which implies that $s\equiv 0\pmod{3}$. There exists a {\it resolvable group divisible $3$-design of type $2^k$} by \cite{BW,brouwer,RS}, namely the $\binom{s}{2}$ edges of a complete graph on $s$ vertices can be partitioned into $\frac{s-2}{2}=k-1$ sets that contain $\frac{s}{3}$ disjoint triangles each, and a $k^{th}$ class, which is a perfect matching. Let $P_1,\ldots,P_{k-1},P_k$ be the partition classes 
where for, $i<k$, $P_i$ consists of the aforementioned disjoint triangles and $P_k$ is the perfect matching. 
Define $H_1,H_2,...,H_{k-1}$ as sets of vertex disjoint $K_3$'s. Further, $H_k$ will consist of the $k=\frac{s}{2}$ matching edges
in $P_k$, with a loop added at both ends of each matching edges.

We will compute the probability that an $ab,cd$ crossing edge pair remains crossed in the $k$-planar drawing, as before. 
\begin{enumerate}
\item Vertices $a,b,c,d$ can map to the vertices of the same matching edge in $2^4$ ways, for $\frac{s}{2}$ edges in $8s$ ways.

\item Edges $ab$ and $cd$ can map to the same edge of a $K_3$ in $(k-1)s\cdot 4$ ways.

\item Edges $ab$ and $cd$ can map to the different edges of a $K_3$ in $(k-1)s\cdot 8$ ways. 
\end{enumerate}
Summing over all possibilities and dividing by $s^4$, the total number of labelings of the four vertices, we obtain that 
$$q = \frac{8}{s^3}+\frac{4k-4}{s^3}+\frac{8k-8}{s^3}=\frac{3k-1}{2k^3}.$$
\end{proof}
\end{subsection}

\begin{subsection}{Proof to Theorem~\ref{th:main2}({\normalfont\ref{th:odd}}): The case when $k$  is odd} \label{odd}

\begin{proof}

Note that for $k=1$, $\frac{2}{k(k+1)}=1$. So we may further assume that $k>1$. We modify our original method by allowing some of the edges of our base graph $K_s^{o}$
to appear in several $H_i$s.

Set $s=k+1$, and note that $s$ is even. It is well known that $K_s$ admits a factorization into  $k$ perfect matchings, $M_1,M_2,...,M_k$. $H_i$ will be obtained from $M_i$ by adding loops to every vertex, so edges that appear in more than one (in fact all) of the $H_i$ are the loops.
We still assign the $\xi(v)$ values 
randomly and uniformly from $\{1,2,\ldots,s\}$ for $v \in V(G)$, but when
an edge of $G$ maps to a loop edge of $K_s^{o}$, we 
randomly and uniformly select an $1\leq i \leq k$ and assign the edge to the $i^{th}$ plane.

In the resulting random $k$-planar drawing of $G$, the probability $q$ with which a crossing edge pair $\{ab,cd\}$ of the optimal planar drawing of $G$ will cross  is still 
independent of the selection of $\{ab,cd\}$ and is an upper bound for $\alpha_k$. 

If $\vert\{\xi(a),\xi(b)\}\vert=\vert\{\xi(c),\xi(d)\}\vert=2$, then the edge-pair remains crossing in the $k$-planar drawing if
$\{\xi(a),\xi(b)\}=\{\xi(c),\xi(d)\}$, which can happen in $\frac{s}{2}\cdot k\cdot 4$ ways, the probability of this is $\frac{2k}{s^3}$.

If $\xi(a)=\xi(b)=\xi(c)=\xi(d)$, then the edge-pair remains crossing in the $k$-planar drawing if
$ab$ and $cd$ are assigned to the same plane. 
The probability of this is $s\cdot \frac{1}{s^4}\cdot \frac{1}{k}=\frac{1}{ks^3}$.

If for some $i\ne j$ we have $i=\xi(a)=\xi(b)$ and $j=\xi(c)=\xi(d)$, then the edge-pair remains crossing in the $k$-planar drawing if
$ab$ and $cd$ is assigned to the $t$-th plane where $\{i,j\}\in M_t$. 
There are $\frac{s}{2}\cdot k$ matching edges, $ab$ and $cd$ can be assigned to different endvertices in 2 ways,
among $s^4$ maps for these 4 vertices, and the images of $ab$ and $cd$ are present in this plane with probability $\frac{1}{k^2}$.
The probability of this case is $\frac{1}{ks^3}$ again.

Finally, if for some $i\ne j$ we have $\{\xi(a),\xi(b)\},\{\xi(c),\xi(d)\}\}=\{\{i,j\},\{i\}\}$ then the edge pair remains crossing if they both get assigned
to the $t$-th plane where $\{i,j\}\in M_t$
There are $sk$
choices for the matching edge $\{i,j\}$ with a distinguished endvertex $i$, 2 ways to choose the edge that maps on the endvertex $i$, 2 ways to 
map the other edge to the matching edge, and the probability that the edge mapped to the loop gets assigned the right plane is $\frac{1}{k}$. 
The probability that the edge pair remains crossing is $\frac{4}{s^3}$.

Summing over all possibilities yields
$\frac{2k}{s^3}+\frac{4}{s^3}+2\cdot \frac{1}{ks^3}=\frac{2}{k(k+1)}.$
\end{proof}
\end{subsection}

\section{Proof to Theorem~\ref{th:mainbip}} \label{bipartite}

Fix a $k>1$ and assume that $G$ is a bipartite graph, with bipartition $A,B$, so that $V(G)=A\cup B$.  

In this section we modify our procedure by changing the base graph $K_s^{o}$ to a complete bipartite graph $K_{k,k}$ with partite sets $\{a_1,\ldots,a_k\}$ and $\{b_1,\ldots,b_k\}$.
 The graphs $H_1,H_2,...,H_k$ are perfect matchings that
make  a factorization of $K_{k,k}$.
 As before, call the components (i.e. edges) of $H_1,H_2,...,H_k$ {\it types}. 

For a vertex $v\in A$, let $\xi(v)$ be a randomly and uniformly distributed value from $\{a_1,...,a_k\}$, and
for a vertex $v\in B$, $\chi(v)$ be a randomly and uniformly distributed value  from $\{b_1,...,b_k\}$.
If $uv$ is an edge of $G$ ($u\in A, v\in B$), then we assign $uv$  edge to the $j$-th plane, if $\{\xi(u),\chi(v)\}\in E(H_j)$. As we factorized a 
complete bipartite graph, there is one and only one such $j$. As before, we draw the $uv$ edges in every plane following the curve representing the $uv$ edge in
an optimal drawing
$\mathcal{D}$ of $G$ in the plane.

Assume that $C=a_\ell b_j$   is an edge of $H_i$. Then clearly 
\[
C^{-1}=\{u\in A: \xi(u)=a_\ell\}\cup \{v\in B: \chi(v)=b_j\} 
\]
is a union of components in $G_i$. We repeat this technique
to reduce the number of crossing edge pairs in the $k$-planar drawing: we translate the subdrawings of $G_i$ on the vertex sets  $\{v\in V(G): \xi(v)\in V(C)\}$
for the components of $H_i$ so far from each other, such that for edges $C_1\not=C_2$ of $H_i$, edges of $G$ between vertices of
$C_1^{-1}$ and edges of $G$ between vertices of $C_2^{-1}$ should not cross.

Assume that $uv$ and $wz$  ($u,w\in A, v,z\in B$) are a pair of crossing edges in the optimal drawing $\mathcal{D}$ of $G$, and have hence 4 distinct
endpoints. The probability that this edge pair is still crossing in the random $k$-planar drawing above, is exactly 
\begin{equation} \label{type1}
q=\PP[{\rm type}(\xi(u),\chi(v))={\rm type}(\xi(w),\chi(z))].
\end{equation}
Note that the value of $q$ does not depend on which crossing edge pair $uv$ and $wz$ was selected from $\mathcal{D}$. Hence 
the expected number of crossings in our random $k$-planar drawing is at most 
$q\Cr(G),$
and therefore some $k$-planar drawing of $G$ has at most  $q\Cr(G)$ crossings. If this holds with a certain $q$ for all bipartite graphs $G$,
then we have established 
\begin{equation} \label{betabound}
\beta_k\leq q.
\end{equation}

If $ab,cd$ are an edge pair that intersects in $\mathcal{D}$, then they remain intersecting in the $k$-planar drawing when they are exactly the same type (\ref{type1}). 
The probability
of that happening is $q=\frac{k^2}{k^4}=\frac{1}{k^2}$, giving $\beta_k\le\frac{1}{k^2}$ by (\ref{betabound}).

Note that the lower bound $\alpha_k\geq 1/k^2$ in Pach et al.~\cite{pach} depends on the existence of the \emph{midrange crossing constant} $\kappa>0$ from Pach, Spencer, and T\'oth~\cite{pachspencertoth}, but not on its value, which is not known. Let $\kappa(n,e)$ denote the minimum crossing number of a
graph $G$ with $n$ vertices and at least $e$ edges. That is,
\begin{equation} \label{minimization}
\kappa(n,e)=\min_{\begin{array}{cc}
|V(G)|=n\\ |E(G)|\ge e \end{array}}\Cr(G).
\end{equation}
Then, according to~\cite{pachspencertoth}, there exists a positive constant $\kappa$,
such that 
the limit
$$\lim_{n\rightarrow\infty}\kappa(n,e)\frac{n^2}{e^3}$$
as $e/n\rightarrow \infty$ and $e=o(n^2)$, 
exists and is equal to $\kappa$.  The existence of such a constant was conjectured by Erd\H os and Guy \cite{EG}.
Czabarka, Reiswig, Sz\'ekely and Wang~\cite{classes} noted, that the existence  of the {midrange crossing constant} for all graphs
can be extended to the existence  of the {midrange crossing constant} $\kappa_{\mathcal C}$  for certain graph classes $\mathcal C$, by requiring $G\in \mathcal{C}$  in (\ref{minimization}),
which may or not be equal to the midrange crossing constant $\kappa$ for all graphs. In fact, 
the best known bounds for $\kappa$ are $0.034\le \kappa\le 0.09$; see \cite{PRTT,ackerman,PT}, while Angelini, Bekos, Kaufmann, Pfister and Ueckerdt
\cite{angelini}  implies that the midrange crossing constant for the class of bipartite graphs is at least $16/289>0.055$, making the  conjecture that these two midrange crossing 
constants differ plausible.

The class of bipartite graphs is such a graph class that admits its midrange crossing constant, and therefore the proof of Pach et al.~\cite{pach} to  $\alpha_k\geq 1/k^2$ 
immediately extends to $\beta_k\geq 1/k^2$.

\section{Theorem~\ref{th:main3}: rectilinear drawings}
We repeat the arguments of Pach et al.~\cite{pach} showing that our new upper bounds apply verbatim to the rectilinear $k$-planar crossing numbers.

The results in this paper on $\alpha_k$ are similarly applicable to $\overline{\alpha}_k$. Specifically, the upper bound starts from a fixed straight-line drawing of $G$ with exactly $\rcro(G)$ crossings. Our randomized procedure decomposes $G$ into $k$ graphs $G_1,\dots , G_{k}$, each of which consists of  vertex-disjoint subgraphs induced by the edge types. 
As the drawings of $G_i$ follow a rectilinear drawing, and translations of drawings of components remain rectilinear and the argument still applies.
 The lower bound relies on the existence of a midrange crossing constant $\overline{\kappa}>0$ for the \emph{rectilinear} crossing number, which is established in~\cite{pachspencertoth} even though the constants $\kappa$ and $\overline{\kappa}$ are not necessarily the same.  Furthermore, our result in Theorem~\ref{th:main} on $\beta_k=1/k^2$ also extends to $\overline{\beta}_k=1/k^2$ and we leave the details to the reader.

\end{document}